\documentclass[10pt]{amsart}

\title{Diameters of Ball Intersections}

\author{Meera Mainkar}
\address{Meera Mainkar\\
Department of Mathematics \\ 
                 Central Michigan University}

\author{Benjamin Schmidt}
\address{Benjamin Schmidt\\
Department of Mathematics \\ 
                 Michigan State University}

\newtheorem{thm}{Theorem}[section]
\newtheorem{lem}[thm]{Lemma}

\theoremstyle{definition}

\newtheorem*{exmp*}{Examples}

\newtheorem*{defn*}{Definition}
\newtheorem*{thm*}{Theorem}	
\newtheorem*{conj*}{Conjecture}

\numberwithin{equation}{section}

\def\Pb{\ifmmode{\Bbb P}\else{$\Bbb P$}\fi}
\def\Z{\ifmmode{\Bbb Z}\else{$\Bbb Z$}\fi}
\def\Q{\ifmmode{\Bbb Q}\else{$\Bbb Q$}\fi}
\def\C{\ifmmode{\Bbb C}\else{$\Bbb C$}\fi}
\def\R{\ifmmode{\Bbb R}\else{$\Bbb R$}\fi}
\def\H{\ifmmode{\Bbb H}\else{$\Bbb H$}\fi}
\def\S{\ifmmode{S^2}\else{$S^2$}\fi}

\def\inj{\operatorname{Inj}}

\def\diam{\operatorname{Diam}}

\def\S{\mathcal S}

\def\conv{\operatorname{Conv}}
\def\conj{\operatorname{Conj}}
\def\foc{\operatorname{Foc}}

\usepackage{graphicx} 

\usepackage{amsmath,amsthm,amsfonts,amscd,flafter,epsf}
\usepackage{graphics}
\usepackage{epsfig}
\usepackage{psfrag}

\begin{document}

\begin{abstract}  
We prove the diameter of the intersection of two closed convex balls in a Riemannian manifold eventually decreases continuously as the centers of the balls move apart.\end{abstract}

\maketitle

\setcounter{secnumdepth}{1}

\section{\bf Prelude}

\begin{center}  Within $R>0$ from Jen and $r>0$ from Jay\\
Poodle and Hound, loyally, do stay.\\
Poodle's utmost insistence?\\
Betwixt Hound, maximal distance.

 \vskip 5pt

Striding subject to such wonder,\\
Weary Jen and Jay asunder. \\
And with their motion set apart,\\
Canines e'er closer, will soon start. 
\end{center}

\setcounter{section}{0}





\section{\bf Statement of Main Theorem}
\textit{As two closed convex balls move apart, the diameter of their intersection eventually decreases continuously}.  This intuitive statement is justified herein.

Let $X$ be a complete Riemannian manifold with positive convexity radius $\conv(X)$ and dimension at least two. The Riemannian structure induces a metric on $X$ $$d:X \times X \rightarrow \mathbb{R}$$ which is complete and geodesic. Denote the closed metric ball with center $x \in X$ and radius $s>0$ by $D^x_s=\{y\,\vert\, d(x,y)\leq s\}.$  Given radii of closed convex balls $0<r\leq R< \conv(X)$ and an arclength parameterized geodesic $\gamma:[0,R+r] \rightarrow X,$ define a function $w:[0,R+r] \rightarrow \mathbb{R}$  by $$w(t)=\diam(D^{\gamma(0)}_R \cap D^{\gamma(t)}_r).$$ The set $$P=\{s \in [0,R+r]\, \vert\, s \leq t \leq R+r \implies D^{\gamma(0)}_R\cap D^{\gamma(t)}_r \subset D^{\gamma(0)}_R\cap D^{\gamma(s)}_r \}$$ is closed and $R+r \in P$.  Therefore $P$ has a well-defined minimum $$T=\min P.$$ 

\noindent \textbf{Main Theorem:} \textit{The restriction of $w$ to $[T,R+r]$ is continuous and strictly decreasing.  Moreover, }

\begin{enumerate} 
\item The restriction of $w$ to $[0,R-r]$ equals $2r$,
\item $R-r\leq T$ with equality if and only if $R=r$,
\item $T<R$,
\item If $R-r<t<R+r$, then $w(t)>R+r-t$. 
\end{enumerate} \vskip 5pt

The convexity assumption is necessary: If $D_1$ and $D_2$ are closed metric balls in the unit sphere $S^2$ of radii at least $\conv(S^2)=\pi/2$, then $\diam(D_1 \cap D_2)=\pi$.

\section{\bf Application of Main Theorem}

Given a self-map $f$ of a metric space $(M,\rho)$, let $$P_f=\{r> 0\, \vert\, \rho(m,\bar{m})=r \implies \rho(f(m),f(\bar{m}))=r\}$$ $$SP_f=\{r > 0\, \vert\, \rho(m,\bar{m})=r \iff \rho(f(m),f(\bar{m}))=r\}.$$  If $(M,\rho)$ is a Euclidean space $\mathbb{E}^d$ with $d\geq 2$, then $P_f =\emptyset$ or $f$ is an isometry \cite{BeQu}.  Other spaces admit non-isometric self-maps with $P_f \neq \emptyset$.\vskip 3pt
\noindent \textbf{Example 1:} Fixing irrationals and shifting rationals one unit right defines a self-map of $\mathbb{E}^1$ with $(0,\infty) \cap \mathbb{Q} \subset SP_f$.\vskip 3pt
\noindent \textbf{Example 2:} Given a subset $A$ of the unit sphere $S^n \subset \mathbb{E}^{n+1}$ with $A=-A$, fixing $A$ and multiplying by $-1$ on the complement of $A$ defines a self-map $f$ of $S^n$ with $\{\frac{1}{2}\pi,\pi\} \subset SP_f$. \vskip 3pt

Conjecturally, dimensional and convexity assumptions exclude self-maps preserving a sufficiently small distance \cite{MaSc}.\vskip 3pt
\noindent\textbf{Conjecture:} \textit{A self-map $f$ of a complete Riemannian manifold $X$ with positive convexity radius $\conv(X)$ and $\dim(X)\geq 2$ satisfies $(0,\conv(X)) \cap P_f =\emptyset$ or is an isometry.}\vskip 3pt

The Conjecture holds for real hyperbolic spaces \cite{Ku} and round spheres \cite{Ev}.  If $f$ is a \textit{bijection} of a locally compact geodesically complete CAT(0) space with path connected metric spheres, then $SP_f=\emptyset$ or $f$ is an isometry \cite{Be, An}; complete and simply connected Riemannian manifolds with nonpositive sectional curvatures are examples of such spaces.  

Additional supporting Theorems are proved in \cite{MaSc}.  In particular, the following is proved using the Main Theorem specialized to the case $r=R$ .\vskip 3pt

\noindent \textbf{Application Theorem:} \textit{For $X$ a connected two-point homogeneous space with $\dim(X) \geq 2$ and $f$ a continuous or surjective self-map of $X$, if $(0,\frac{2}{3}\conv(X)) \cap SP_f \neq \emptyset$, then $f$ is an isometry.}\vskip 3pt

The connected two-point homogenous spaces consist precisely of the Euclidean spaces $\mathbb{E}^n$ and the rank one symmetric spaces $\mathbb{R}H^n$, $\mathbb{C}H^n$, $\mathbb{H}H^n$, $\mathbb{O}H^2$, $S^n$, $\mathbb{R}P^n$, $\mathbb{C}P^n$, $\mathbb{H}P^n$, and $\mathbb{O}P^2$ \cite{Wa,Sz}.  The unified proof does not use this classification.

  \section{\bf Tools} Results used in the proof of the Main Theorem are now summarized.  Possible references include \cite{doCa, BuBuIv, Di}.

\subsection{Notation}
Given a metric space $(M,\rho)$, a nonempty closed subset $Y\subset M$, and $s\in (0,\infty)$,  let $$B^Y_s=\{m\,\vert\,\rho(m,Y)<s\},\hskip 5pt D^Y_s=\{m\,\vert\, \rho(m,Y)\leq s\},\hskip 5pt S^Y_s=\{m\,\vert\,\rho(m,Y)=s\}.$$ When $Y=\{y\}$, braces $\{\}$ are omitted in the simpler notation $B^{y}_s$, $D^y_s$, and $S^y_s$.

\subsection{Hausdorff Distance} 
Given nonempty closed subsets $Y, Z \subset M$, the \textit{Hausdorff distance} between $Y$ and $Z$, denoted $\rho_H(Y,Z)$, is defined by $$\rho_H(Y,Z)=\inf\{s>0\,\vert\, Y \subset B^Z_s\,\,\text{and}\,\,Z \subset B^Y_s\}.$$ The Hausdorff distance defines a metric on the set of nonempty, closed, and bounded subsets of $M$.  The next two Lemmas are known and easily proved.

\begin{lem}\label{diameter}
If $Y$ and $Z$ are compact, then $|\diam(Y)-\diam(Z)|\leq 2\rho_H(Y,Z).$
\end{lem}

\begin{lem}\label{converge}
Let $\{Y_i\}_{i=1}^{\infty}$ be a sequence of compact subsets.
\begin{enumerate}
\item If $Y_{i+1} \subset Y_i$ for each $i$, then $\displaystyle\lim_{k\rightarrow \infty} Y_k$$=\bigcap_i Y_i$ 
\item If $Y_i \subset Y_{i+1}$ for each $i$ and $\bigcup_i Y_i$ is compact, then $\displaystyle\lim_{k\rightarrow \infty}Y_k$$=\bigcup_i Y_i$. 
\end{enumerate}
\end{lem}

\subsection{Riemannian Distance and Geodesic Variations}Let $(X,g)$ denote a complete Riemannian manifold.  For $I \subset \mathbb{R}$ an interval, let $|I|$ denote its length.  For $p \in X$, let $T_pX$ denote the tangent space to $X$ at $p$.  If $v\in T_pX$, let $||v||=g_p(v,v)^{\frac{1}{2}}$ denote its length.  Given a (piecewise) smooth curve $$c:I\rightarrow X,$$ its energy $E(c)$ and length $L(c)$ are defined by $$E(c)=\int_I ||\dot{c}(t)||^2 dt\,\,\,\,\,\,\,\,\,\,\,\,\,\,\text{and}\,\,\,\,\,\,\,\,\,\,\,\,\, L(c)=\int_I ||\dot{c}(t)||\,dt.$$  The Cauchy-Schwartz inequality implies that $$L^2(c) \leq |I| E(c)$$  with equality holding if and only if the speed  $||\dot{c}(t))||$ is constant.  The length functional on paths equips $X$ with a complete metric $$d:X \times X \rightarrow \mathbb{R}$$ for which $d(p,q)$ equals the minimal length of a smooth path in $X$ starting at $p$ and ending at $q$.  A \textit{geodesic} is a smooth, constant speed, path $c:I \rightarrow X$ satisfying $$L(c_{\vert [t_1,t_2]})=d(c(t_1),c(t_2))$$ for all subintervals $[t_1,t_2]\subset I$ of sufficiently short length.  A geodesic is \textit{minimizing} if the former holds for all subintervals of $I$.  In particular, if  $c:[0,1] \rightarrow X$ is a \textit{minimizing geodesic},  then \begin{equation} \label{energy} d^2(c(0),c(1))=E(c).\end{equation}

If $a,b,c$ are three points in $X$ that do not lie in the image of a common minimizing geodesic then the following \textit{strict triangle inequality} holds: \begin{equation}\label{strict}d(a,b) < d(a,c)+d(c,b). \end{equation}

Fix $p \in X$.  Each $v \in T_pX$ determines a unique geodesic $c_v:\mathbb{R} \rightarrow X$ with $\dot{c}_v(0)=v$.  The exponential map $\exp_p:T_pX \rightarrow X$ is defined by $\exp_p(v)=c_v(1)$.  Let $\epsilon>0$ and let $$V:(-\epsilon,\epsilon) \rightarrow T_pX$$ be a smooth path. Consider the parameterized surface $$f:[0,1]\times (-\epsilon,\epsilon) \rightarrow X$$ defined by $$f(t,s)=\exp_p(tV(s)).$$  Consider the curves $$\sigma(s)=f(1,s),\,\,\,\,\,c_s(t)=f(t,s),\,\,\,\,\,\text{and}\,\,\,\,\, c(t)=c_0(t),$$ and the vector fields $$f_t=df(\frac{\partial}{\partial t}),\,\,\,\,\, f_s=df(\frac{\partial}{\partial s}),\,\,\,\,\,\text{and}\,\,\,\,\, J(t)=f_s(c(t)).$$  Then for each $s \in (-\epsilon,\epsilon)$, $$c_s:[0,1] \rightarrow X$$ is a geodesic.  The vector field $J(t)$ satisfies $J(0)=0$ and the following defining equality of a \textit{Jacobi field along the geodesic} $c(t)$:  $$J''+R(J,\dot{c})\dot{c}=0.$$ If additionally, $g(J,\dot{c})=0,$ then $J(t)$ is said to be a \textit{normal Jacobi field along the geodesic} $c(t)$. The first and second derivative formulas for $s \mapsto E(c_s)$ at $s=0$ are given by:

\begin{equation}\label{first}
\frac{\mathrm{d}E(c_s)}{\mathrm{d}s}(0)=2g(\dot{\sigma}(0),\dot{c}(1))
\end{equation}

\begin{equation}\label{second}
\frac{\mathrm{d}^2E(c_s)}{\mathrm{d} s^2}(0)=\frac{\mathrm{d} ||J||^2}{\mathrm{d} t}(1)+2g(\nabla_{J(1)} f_s,\dot{c}(1)).
\end{equation}

\subsection{Riemannian Convexity}
As above, $(X,g)$ denotes a complete Riemannian manifold. Let $I \subset \mathbb{R}$ be an interval and $\mathcal{O} \subset X$ an open subset.  A real valued function $$F: \mathcal{O} \rightarrow \mathbb{R}$$ is \textit{strictly convex} if for each non-constant geodesic $$\tau:I \rightarrow \mathcal{O},$$ the function $$h:I \rightarrow \mathbb{R}$$ defined by $h=F\circ \tau$ is \textit{strictly convex}: For each distinct $s,t \in I$ and $\lambda \in (0,1)$, $$h(\lambda s + (1-\lambda)t)<\lambda h(s)+(1-\lambda)h(t).$$  When $F$ is twice continuously differentiable, the latter is equivalent to the second derivative inequality $h''<0$.

A subset $Y \subset X$ is \textit{strongly convex} provided that whenever $y_1,y_2 \in Y$, there exists a unique minimizing geodesic in $X$ with endpoints $y_1$ and $y_2$, and moreover, this geodesic lies entirely in $Y$. Following the presentation in \cite{Di}, we now define several metric invariants of $X$.

The \textit{convexity radius of a point} $x\in X$, denoted $\conv(x)$, is the supremum of positive real numbers $s>0$ having the property that for each $0<r<s$, the open ball $B^x_r$ is strongly convex.  The \textit{convexity radius of} $X$, denoted $\conv(X)$, is the infimum of the convexity radii of its points.  

The \textit{injectivity radius of a point} $x\in X$, denoted $\inj(x)$, is the supremum of positive real numbers $s>0$ having the property that the restriction of the exponential map $\exp_x:T_x X \rightarrow X$ to the open ball of radius $s$ and center $0$ is a diffeomorphism onto its image $B^x_s$.   The \textit{injectivity radius of} $X$, denoted $\inj(X)$, is the infimum of the injectivity radii of its points.  Every geodesic starting at $x \in X$ of length less than $\inj(x)$ is minimizing.

The \textit{conjugate radius of a point} $x \in X$, denoted $\conj(x)$, is the minimum $T>0$ such that there exists a unit-speed geodesic $c:\mathbb{R} \rightarrow X$ and a non-zero normal Jacobi field $J(t)$ along $c(t)$ with $c(0)=x,$ $J(0)=0$, and $J(T)=0$.  If no such $T$ exists, then the conjugate radius at $x$ is infinite.  The \textit{conjugate radius of} $X$, denoted $\conj(X)$, is the infimum of the conjugate radii of its points.

The \textit{focal radius of a point} $x\in X$, denoted $\foc(x)$, is the minimum $T>0$ such that there exists a unit-speed geodesic $c:\mathbb{R} \rightarrow X$ and a non-zero normal Jacobi field $J(t)$ along $c(t)$ with $c(0)=x$, $J(0)=0$, and $\frac{\mathrm{d}||J||}{\mathrm{d}t}(T)=0.$ If no such $T$ exists, then the focal radius at $x$ is infinite.  The \textit{focal radius of} $X$, denoted $\foc(X)$, is the infimum of the focal radius of its points.

A \textit{geodesic loop in} $X$ based at a point $x \in X$ is a geodesic $c:[0,1] \rightarrow X$ with $c(0)=x=c(1)$.  Given $x \in X$, let $L(x)$ denote the minimal length of a non-constant geodesic loop based at $x \in X$.  If no such loop exists, then $L(x)$ is infinite.  Let $L(X)$ denote the infimum of $L(x)$ over points $x \in X$.

\begin{lem}\label{convex}
Let $X$ be a complete Riemannian manifold with positive convexity radius $\conv(X)$ and with induced metric $d:X \times X \rightarrow \mathbb{R}$. Let $l \in (0,\conv(X))$ and $u \in X$.  Then
\begin{enumerate}
 \item $\conv(X)=\min\{\foc(X), \inj(X)/2\}$,\vskip 3pt
 \item $d^2(u,\cdot): B^u_{\conv(X)} \rightarrow [0,\conv(X)^2)$ is strictly convex,\vskip 3pt
 \item If $\sigma:(-\epsilon,\epsilon) \rightarrow X$ is a non-constant geodesic with $\dot{\sigma}(0)$ tangent to the sphere $S^{u}_l$, then for all non-zero $s$ sufficiently close to $0$, $d(u, \sigma(s))>l$, \vskip 3pt
 \item If $p, q\in D^{u}_l$ satisfy $0<d(p,q)<2l$, then for each $\alpha>0$ there exists $a,b \in D^u_l$ with  $a\in B^p_{\alpha}$,  $b \in B^q_{\alpha}$, and $d(a,b)>d(p,q)$, \vskip 3pt
 \item If $c:[0,1] \rightarrow X$ is a non-constant geodesic with $c(0) \in D^{u}_l$ and $c(1) \in S^{u}_l$ and $\gamma:\mathbb{R} \rightarrow X$ is the geodesic with $\gamma(0)=u$ and $\gamma(l)=c(1)$, then $g(\dot{c}(1),\dot{\gamma}(l))>0$, and\vskip 3pt
 \item The closed metric ball $D^{u}_l$ is strongly convex.
  \end{enumerate}
 \end{lem}
 
 \noindent{\textit{Proof of (1).}} By \cite{Kl} $$\frac{\inj(X)}{2}=\min\{\frac{\conj(X)}{2},\frac{L(X)}{4}\}.$$  By \cite{Di} $$\conv(X)=\min\{\foc(X),\frac{L(X)}{4} \} .$$  It remains to observe that $$\foc(X) \leq \frac{\conj(X)}{2}.$$ This is immediate since if $c:\mathbb{R} \rightarrow X$ is a unit-speed geodesic, $S>0$, and $J(t)$ is a non-zero normal Jacobi field along $c(t)$ with $J(0)=0=J(S)$, then $||J||(t)$ is maximized at some parameter $\bar{S} \in (0,S)$ and $$\min\{\bar{S},S-\bar{S}\}\leq \frac{S}{2}.$$ \qed \vskip 5pt

 \noindent{\textit{Proof of (2).}} The proof is an application of the second derivative formula (\ref{second}) as stated in \cite{Di}; we include a proof .  Let $$\sigma:(-\epsilon,\epsilon)\rightarrow B^{u}_{\conv(X)}$$ be a non-constant geodesic.  By Item (1), $\conv(X)< \inj(X)$ so that  the restriction of $\exp_u$ to $B^{0}_{\conv(X)}$ is a diffeomorphism onto its image $B^{u}_{\conv(X)}$. Define $$V:(-\epsilon,\epsilon) \rightarrow B^0_{\conv(X)}\subset T_uX$$ by $$\sigma(s)=\exp_u(V(s))$$ and consider the smooth geodesic variation $$f:[0,1]\times (-\epsilon,\epsilon) \rightarrow X$$ defined by $$f(t,s)=\exp_u(tV(s)).$$  Let $c_s(t)=f(t,s)$, $c(t)=c_0(t)$, and let $J(t)$ be the associated Jacobi field along $c(t)$.  As $\conv(X)<\inj(X)$, each geodesic $c_s$ is minimizing.  By (\ref{energy}) $$d^2(u,\sigma(s))=E(c_s).$$ Note that $$\nabla_{J(1)}f_s=\nabla_{\dot{\sigma}} {\dot{\sigma}}(0)$$ and since $\sigma$ is a geodesic, the second term in (\ref{second}) is zero. Hence, $$\frac{\mathrm{d}^2E(c_s)}{\mathrm{d} s^2}(0)=\frac{\mathrm{d} ||J||^2}{\mathrm{d} t}(1)=2||J||(1) \cdot \frac{\mathrm{d} ||J||}{\mathrm{d} t}(1).$$ The latter is positive since by Item (1), $\conv(X)\leq \foc(X)$, concluding the proof. \qed \vskip 5pt
 
 \noindent{\textit{Proof of (3).}} By (1), $l<\inj(X)$.  Therefore $S^{u}_l$ is a smoothly embedded submanifold of $X$ so that the statement is meaningful. After reducing $\epsilon$ if necessary, Item (2) implies the function $h(s)=d^{2}(u,\sigma(s))$ is strictly convex and satisfies $h(0)=l^2$.  It therefore suffices to prove that $h'(0)=0$.
 
For each $s \in (-\epsilon,\epsilon)$ there is a unique minimizing geodesic $$c_s:[0,1]\rightarrow X$$ joining $c_s(0)=u$ to $c_s(1)=\sigma(s)$. By (\ref{energy}) and (\ref{first}), $h'(0)=2g(\dot{c}_0(1),\dot{\sigma}(0))$.  The latter is zero since by Gauss' Lemma $\dot{c}_0(1)$ is perpendicular to $T_{\sigma(0)} S^{u}_l.$ \qed \vskip 5pt

\noindent{\textit{Proof of (4).}} Let $c:\mathbb{R} \rightarrow X$ be the complete geodesic whose restriction to $[0,1]$ is a minimizing geodesic joining $c(0)=p$ to $c(1)=q$.    The hypothesis and Item (1) imply $$d(p,q)<2l<\inj(X).$$

We first consider the case when $q \in B^u_l$.  For all $t_0$ sufficiently close to $1$, $$c(t_0) \in B^u_l\cap B^{q}_{\alpha}$$ and $$d(p,c(t_0))<\inj(X).$$ Choosing such a $t_0$ that is also greater than $1$,  $$d(p,c(t_0))=L(c_{\vert [0,t_0]})>L(c_{\vert [0,1]})=d(p,q).$$ Setting $p=a$ and $b=c(t_0)$ completes the proof in this case.

 We now consider the case when $q \notin B^{u}_l$.  In this case, $q \in S^{u}_l$.  The metric sphere $S^{u}_l$ is a smooth codimension one submanifold of $X$ since by Item (1), $l<\inj(X)$.  In particular, it has a well-defined tangent space at each point.

 Consider the case when $\dot{c}(1)$ is perpendicular to $T_q S^{u}_l$. By (1), there exists a unique length $l$ vector $v\in T_uX$ with $$\exp_u(v)=c(1).$$  Let $$c_v:\mathbb{R} \rightarrow X$$ denote the complete geodesic with $\dot{c}_v(0)=v$. Gauss' Lemma implies that the image of $c:[0,1] \rightarrow X$ is a subset of the image of $c_v:[-1,1] \rightarrow D^{x}_r$, a geodesic segment of length $2l$ (with midpoint $u$).  As $d(p,q)<2l$, it follows that $$p\in B^{u}_l.$$  The previous argument then applies to complete the proof in this case.

Finally, if $\dot{c}(1)$ isn't perpendicular to $T_q S^{u}_l$, then there exists $v \in T_q S^{u}_l$ with $$g(\dot{c}(1),v)>0.$$  Choose a smooth curve $$\sigma:(-\epsilon,\epsilon) \rightarrow S^u_l$$ with $$(\sigma(0),\dot{\sigma}(0))=(q,v)$$ and with $$d(p,\sigma(s))<\inj(X)$$ for each $s \in (-\epsilon,\epsilon)$.  Then for each $s \in (-\epsilon,\epsilon)$, there exists a unique minimizing geodesic $$c_s:[0,1] \rightarrow X$$ joining $c_s(0)=p$ to $c_s(1)=\sigma(s)$.  Note that $c_0$ is the restriction of $c$ to $[0,1]$.  By (\ref{energy}) and (\ref{first}), if $s_0$ is positive and sufficiently close to zero, then $$d(p,\sigma(s_0))>d(p,q)\,\,\,\,\, \text{and}\,\,\,\,\, \sigma(s_0)\in B^q_{\alpha} .$$  Setting $p=a$ and $b=\sigma(s_0)$ completes the proof. \qed \vskip 5pt

\noindent{\textit{Proof of (5).}}  By (2), $f(t)=d^2(u,c(t))$ is strictly convex on $[0,1]$.  As $f(0)\leq f(1)$, $f'(1)>0$.  The desired inequality now follows from (\ref{first}).\qed \vskip 5pt

\noindent{\textit{Proof of (6).}}  Let $p, q \in D^{u}_l$.  As $D^{u}_l \subset B^{u}_{\conv(X)}$ there is a unique minimizing geodesic $\gamma:[0,1] \rightarrow X$ with $\gamma(0)=p$ and $\gamma(1)=q$.  Moreover $\gamma$ has image in $B^{u}_{\conv(X)}.$  By (2), $f(t)=d^2(u,\gamma(t))$ is strictly convex.  Therefore, for each $t \in (0,1)$, $d^2(u,\gamma(t))<\max\{d^2(u,p),d^2(u,q)\}\leq l^2,$ demonstrating that $\gamma$ has image in $D^{u}_l$.\qed

\section{\bf Proof of Main Theorem}  By Lemma \ref{convex}-(1), the arclength parameterized geodesic $\gamma:[0,R+r] \rightarrow X$ is minimizing: For each $t_1,t_2 \in [0,R+r]$, 
\begin{equation}\label{geo}
d(\gamma(t_1),\gamma(t_2))=|t_1-t_2|.
\end{equation} \vskip 5pt

Items (1)-(4) in the Main Theorem are proved first.

\subsection{Proof of Item (1): The restriction of $w$ to $[0,R-r]$ equals $2r$.\\ \\}  
Let $t \in [0,R-r]$ and $x \in D^{\gamma(t)}_r$.  By the triangle inequality and (\ref{geo}) $$d(\gamma(0),x)\leq d(\gamma(0), \gamma(t))+d(\gamma(t),x) \leq t+r \leq R.$$ Therefore $$D^{\gamma(t)}_r \subset D^{\gamma(0)}_R$$ and $$w(t)=\diam(D^{\gamma(0)}_R\cap D^{\gamma(t)}_r)=\diam(D^{\gamma(t)}_r)=2r.$$ \qed \vskip 5pt

\subsection{Proof of Item (2): The inequality $R-r\leq T$ holds and is an equality if and only if $R=r$.\\ \\} 
When $R=r$, it is trivial to verify that $0 \in P$, whence $$T=0=R-r.$$  

Now assume that $R>r$.  To demonstrate that $T>R-r$, it suffices to show that for each $s \in [0,R-r]$, there exists $x \in D^{\gamma(0)}_R \cap D^{\gamma(R)}_r$ satisfying $r<d(\gamma(s),x).$ If $s=0$, let $x=\gamma(R)$.  By (\ref{geo}), $$r<R=d(\gamma(0),\gamma(R))=d(\gamma(0),x),$$ concluding the proof in this case.

If $s>0$, then choose $x \in S^{\gamma(0)}_R \cap D^{\gamma(R)}_r$ distinct from $\gamma(R)$.  Note that $x$ is distinct from $\gamma(-R)$ since $$d(x,\gamma(R))\leq r<2R=d(\gamma(-R),\gamma(R)).$$  In particular, the strict triangle inequality (\ref{strict}) applies to the triple $(\gamma(0),\gamma(s),x)$.  The strict triangle inequality and (\ref{geo}) imply $$R=d(\gamma(0),x)< d(\gamma(0),\gamma(s))+d(\gamma(s),x)= s+d(\gamma(s),x).$$ Therefore $$d(\gamma(s),x)>R-s\geq r,$$  concluding the proof. \qed \vskip 5pt

\subsection{Proof of Item (3): The inequality $T<R$ holds.\\ \\} 
The proof is based on the following Claim. \vskip 5pt

\noindent \underline{Claim:} \textit{There exists $v\in T_{\gamma(T)}X$ with $g(v,\dot{\gamma}(T)) \geq 0$, $\exp_{\gamma(T)}(v) \in D^{\gamma(0)}_R$, and $||v||=r$.}\vskip 5pt

\noindent \underline{Proof of Item (3) assuming Claim:}  We first argue that $R \in P$.  To this end, let $R<t\leq R+r$ and let $$x \in D^{\gamma(0)}_R\cap D^{\gamma(t)}_r.$$ We must demonstrate $d(x,\gamma(R)) \leq r$. As $$\gamma(0), \gamma(t) \in B^{x}_{\conv(X)},$$ a strongly convex ball, the restriction of $\gamma$ to $[0,t]$ has image in $B^{x}_{\conv(X)}$.  By Lemma \ref{convex}-(2), the function $$f(\cdot)=d^2(x,\gamma(\cdot))$$ is strictly convex on $[0,t]$.  If $$c:[0,1] \rightarrow X$$ is a minimizing geodesic joining $c(0)=x$ to $c(1)=\gamma(R)$, then by Lemma \ref{convex}-(5), $$g(\dot{c}(1),\dot{\gamma}(R))>0.$$ Therefore, by (\ref{first}), $f'(R)>0$.  As $f$ is convex, $f$ is increasing on $[R,t]$ whence $$d^2(x,\gamma(R))=f(R)<f(t)=d^2(x,\gamma(t))\leq r^2,$$ concluding the proof that $R\in P$. 

As $R \in P$, $T=\min{P}\leq R$.  We conclude by showing that $T \neq R$.  According to the Claim, if $T=R$, then there exists $v\in T_{\gamma(R)}X$ with $||v||=r$, $g(v, \dot{\gamma}(R)) \geq 0$, and $\exp_{\gamma(R)}(v) \in D^{\gamma(0)}_R.$  As $$\{\gamma(R), \exp_{\gamma(R)}(v)\} \subset D^{\gamma(0)}_R,$$ Lemma \ref{convex}-(6) implies that for each $t \in [0,1]$, $$\exp_{\gamma(R)}(tv) \in D^{\gamma(0)}_R.$$  A contradiction is obtained when $g(v, \dot{\gamma}(R))>0$ since the gradient of $d(\gamma(0),\cdot)$ at $\gamma(R)$ equals $\dot{\gamma}(R)$ and when $g(v,\dot{\gamma}(R))=0$ by Lemma \ref{convex}-(3). \qed \vskip 5pt

\noindent \underline{Proof of Claim.} If $T=0$, then any vector of length $r$ perpendicular to $\dot{\gamma}(0)$ satisfies the conditions in the Claim.  

Now suppose that $T>0$.  For each $n\in \mathbb{N}$ satisfying $\frac{1}{n}<T$, let $$s_n=T-\frac{1}{n}.$$  By definition of $T$, there exists $$\bar{t}_n>s_n\,\,\,\,\, \text{and}\,\,\,\,\, x_n \in D^{\gamma(0)}_R$$ such that \begin{equation}\label{blah} d(x_n,\gamma(\bar{t}_n)) \leq r<d(x_n,\gamma(s_n)).\end{equation}

For each index $n$, define $t_n$ by $t_n=\bar{t}_n$ if $\bar{t}_n<T$ and by $t_n=T$ if $\bar{t}_n \geq T$.  Then for each index $n$,  (\ref{geo}) implies \begin{equation} \label{s} d(\gamma(t_n),\gamma(s_n))=t_n-s_n\leq \frac{1}{n}.\end{equation}  Moreover, for each index $n$, \begin{equation}\label{sh} d(x_n,\gamma(t_n)) \leq r<d(x_n,\gamma(s_n))\end{equation} by (\ref{blah}) and the fact that $T \in P$.  The triangle inequality and (\ref{s}-\ref{sh}) now imply \begin{equation}\label{close}  r<d(x_n,\gamma(s_n))\leq d(x_n,\gamma(t_n))+d(\gamma(t_n),\gamma(s_n)) \leq r+\frac{1}{n} \end{equation} Let $$c_n:[0,1] \rightarrow X$$ be a minimizing geodesic with $$c_n(0)=\gamma(s_n)\,\,\,\,\, \text{and}\,\,\,\,\, c_n(1)=x_n.$$  Let $v$ be a limit point of the sequence of vectors $\{\dot{c_n}(0)\}$.  Then $v$ is tangent to $X$ at $\gamma(T)$ and by (\ref{close}), $||v||=r$.  As $D^{\gamma(0)}_R$ is closed and $x_n \in D^{\gamma(0)}_R$ for each $n$, it follows $$\exp_{\gamma(T)}(v) \in D^{\gamma(0)}_R.$$  

To prove the remaining assertion that $g(v,\dot{\gamma}(T)) \geq 0$ it suffices to prove that $g(\dot{c}_n(0),\dot{\gamma}(s_n) \geq 0$ for sufficiently large $n$. For $n$ sufficiently large, $$r+\frac{1}{n}<\conv(X).$$ By (\ref{close}) and Lemma \ref{convex}, the function $$f_n:[s_n,t_n] \rightarrow  \mathbb{R}$$ defined by $$f_n(t)=d^2(x_n,\gamma(t))$$ is strictly convex.  By (\ref{sh}), $f_n$ is initially decreasing.  The desired inequality now follows from (\ref{first}). \qed \vskip 5pt

\subsection{Proof of Item (4): If $R-r<t<R+r$, then $w(t)>R+r-t$.\\ \\} 

Fix $t \in (R-r,R+r)$ and set $$T_0=\frac{R-r+t}{2}\,\,\,\,\, \text{and}\,\,\,\,\, T_1=R-T_0=\frac{R+r-t}{2}.$$  Let $c:\mathbb{R} \rightarrow X$ denote the complete geodesic uniquely determined by $$(c(0),\dot{c}(0))=(\gamma(T_0),\dot{\gamma}(T_0)).$$


The inequalities $0\leq R-r<t<R+r$  and (\ref{geo}) imply that $$c(T_1)=\gamma(T_0+T_1)=\gamma(R) \in S^{\gamma(0)}_R \cap B^{\gamma(t)}_r$$ and $$c(-T_1)=\gamma(T_0-T_1)=\gamma(t-r) \in B^{\gamma(0)}_R \cap S^{\gamma(t)}_r. $$  As $\gamma$ is parameterized by arclength, so too is $c$.  Therefore, the restriction of $c$ to $[-T_1,T_1]$ has length $$2T_1=R+r-t.$$  By Lemma \ref{convex}-(1) $$d(c(-T_1),c(T_1))=R+r-t.$$  Conclude that $$w(t) \geq R+r-t.$$  Next, we demonstrate that this inequality is strict.

Let $v$ be a unit length tangent vector to $X$ at $c(0)=\gamma(T_0)$ which is close to but distinct from $\dot{c}(0)=\dot{\gamma}(T_0).$  Set $\tau(s)=\exp_{c(0)}(sv)$.   There exist $T_{+}$ and $T_{-}$ close to $T_1$ and $-T_1$, respectively, such that $$\tau(T_+) \in S^{\gamma(0)}_R \cap B^{\gamma(t)}_r$$ and $$\tau(T_{-}) \in B^{\gamma(0)}_R \cap S^{\gamma(t)}_r. $$  As above, conclude that if $v$ is sufficiently close to $\dot{c}(0)=\dot{\gamma}(T_0)$, then $$w(t) \geq d(\tau(T_+),\tau(T_{-}))=T_{+}-T_{-}.$$  We conclude by arguing that $$T_1<\min\{-T_{-},T_{+}\}.$$   As $v$ is distinct from $\dot{\gamma}(T_0)$, the strict triangle inequality (\ref{strict}) implies that $$d(\gamma(0), \tau(T_{+}))<d(\gamma(0),c(0))+d(c(0),\tau(T_{+}))=d(\gamma(0),\gamma(T_0))+T_{+}.$$ Therefore, $$T_1=R-T_0<T_{+}.$$  Similarly, the strict triangle inequality implies that $$r=d(\tau(T_{-}),\gamma(t))<d(\tau(T_{-}),\tau(0))+d(\tau(0),\gamma(t))=-T_{-}+(t-T_0).$$ Therefore $$T_1=r+T_0-t<-T_{-}.$$  \qed


\subsection{Remainder of Proof: The restriction of $w$ to $[T,R+r]$ is continuous and strictly decreasing.\\ \\} 
These facts are easily deduced from Claims 1-3 below. \vskip 5pt

\noindent \underline{Claim 1:} \textit{The restriction of $w$ to $(T,R+r]$ is left continuous.}\vskip 5pt

\noindent \underline{Claim 2:} \textit{The restriction of $w$ to $[T,R+r)$ is right continuous.}\vskip 5pt

\noindent \underline{Claim 3:} \textit{The restriction of $w$ to $(T,R+r]$ is strictly decreasing.}\vskip 5pt

The proofs of Claims 1-3 use the following Claim. \vskip 5pt

\noindent \underline{Claim:} \textit{If $T\leq s <t \leq R+r$, then $D^{\gamma(0)}_R\cap D^{\gamma(t)}_r \subset D^{\gamma(0)}_R\cap B^{\gamma(s)}_r\subset D^{\gamma(0)}_R\cap D^{\gamma(s)}_r$.}\vskip 5pt

\noindent \underline{Proof of Claim.}
Let $x\in D^{\gamma(0)}_R\cap D^{\gamma(t)}_r$.  We must show that $d(x,\gamma(s))<r$.  By assumption \begin{equation}\label{assume}d(x,\gamma(t))\leq r\end{equation} and by definition of $T$, \begin{equation}\label{df} d(x,\gamma(T))\leq r.\end{equation} By (\ref{assume})-(\ref{df}), $$\gamma(T),\gamma(t)\in B^{x}_{\conv(X)}.$$  As this ball is strongly convex, the restriction of $\gamma$ to $[T,t]$ has image in $B^{x}_{\conv(X)}$.  By Lemma \ref{convex}-(2),  $$d^2(x,\cdot):B^{x}_{\conv(X)} \rightarrow \mathbb{R}$$ is strictly convex. Therefore, $$d^2(x,\gamma(s))<\max\{d^2(x,\gamma(T)),d^2(x,\gamma(t)\}\leq r^2.$$ \qed \vskip 5pt

\noindent \underline{Proof of Claim 1.}
Fix $t \in (T,R+r]$ and $\epsilon>0$.  Lemma \ref{converge}-(1) implies that there exists $$0<\delta<t-T$$ such that 
\begin{equation}\label{estimate2}
d_H(D^{\gamma(0)}_R\cap D^{\gamma(t)}_r, D^{\gamma(0)}_R\cap D^{\gamma(t)}_{r+\delta})<\epsilon/2.
\end{equation} Let $s \in (t-\delta,t)$.  By the Claim, \begin{equation}\label{contain3} D^{\gamma(0)}_R \cap D^{\gamma(t)}_r \subset D^{\gamma(0)}_R \cap D^{\gamma(s)}_r. \end{equation} The triangle inequality and (\ref{geo}) imply \begin{equation}\label{contain4} D^{\gamma(0)}_R \cap D^{\gamma(s)}_{r} \subset D^{\gamma(0)}_R \cap D^{\gamma(t)}_{r+\delta}.\end{equation}
Equations (\ref{estimate2})-(\ref{contain4}) imply $$d_H(D^{\gamma(0)}_R \cap D^{\gamma(s)}_r,D^{\gamma(0)}_R \cap D^{\gamma(t)}_r)<\epsilon/2.$$ Lemma \ref{diameter} implies $$|w(s)-w(t)|<\epsilon.$$ \qed \vskip 5pt

\noindent \underline{Proof of Claim 2.}
Fix $s \in [T,R+r)$ and let $\epsilon>0$.  Lemma \ref{converge}-(2) implies that there exists $$0< \delta<\min\{r,R+r-s\}$$ such that 
\begin{equation}\label{estimate1} d_H(D^{\gamma(0)}_R\cap D^{\gamma(s)}_r, D^{\gamma(0)}_R\cap D^{\gamma(s)}_{r-\delta})<\epsilon/2.\end{equation}  Let $t \in (s,s+\delta)$.  By the Claim, \begin{equation}\label{contain1} D^{\gamma(0)}_R \cap D^{\gamma(t)}_r \subset D^{\gamma(0)}_R \cap D^{\gamma(s)}_r. \end{equation} The triangle inequality and (\ref{geo}) imply \begin{equation}\label{contain2} D^{\gamma(0)}_R \cap D^{\gamma(s)}_{r-\delta} \subset D^{\gamma(0)}_R \cap D^{\gamma(t)}_{r}.\end{equation}
Equations (\ref{estimate1})-(\ref{contain2}) imply $$d_H(D^{\gamma(0)}_R \cap D^{\gamma(s)}_r,D^{\gamma(0)}_R \cap D^{\gamma(t)}_r)<\epsilon/2.$$ Lemma \ref{diameter} implies $$|w(s)-w(t)|<\epsilon.$$ \qed \vskip 5pt

\noindent \underline{Proof of Claim 3.} Assume that $T < s<t \leq R+r$.  Choose $p,q \in D^{\gamma(0)}_R\cap D^{\gamma(t)}_r$ with $$w(t)=d(p,q).$$ By the Claim, $p,q \in D^{\gamma(0)}_R \cap B^{\gamma(s)}_r$.  In particular, $$d(p,q)\leq d(p,\gamma(s))+d(\gamma(s),q)<2r\leq2R.$$

Choose $\alpha>0$ such that $\alpha<\min\{r-d(p,\gamma(s)),r-d(q,\gamma(s))\}$.  By Lemma \ref{convex}-(4), there exists $$a \in D^{\gamma(0)}_R \cap B^{p}_\alpha\,\,\,\,\, \text{and}\,\,\,\,\, b\in D^{\gamma(0)}_R\cap B^{q}_{\alpha}$$ such that $$d(a,b)>d(p,q).$$ By the triangle inequality, $$a \in D^{\gamma(0)}_R\cap D^{\gamma(s)}_r\,\,\,\,\, \text{and}\,\,\,\,\,b \in D^{\gamma(0)}_R\cap D^{\gamma(s)}_r.$$   Therefore, $$w(s)=\diam(D^{\gamma(0)}_R \cap D^{\gamma(s)}_r) \geq d(a,b)>d(p,q)=w(t).$$ \qed \vskip 5pt




\section{\bf Speculation}

Note $w(t)=\diam(D^{\gamma(0)}_R\cap D^{\gamma(t)}_r)\leq \diam(D^{\gamma(t)}_r)=2r$.  Set $$Q=\{t \in [0,R+r]\,\vert\, w(t)=2r\}.$$  The set $Q$ is closed and $0 \in Q$.  Set $S=\max Q$ and note $R-r\leq S\leq T$.

When $X$ has constant sectional curvatures, $S=T$, the function $w$ equals $2r$ on $[0,T]$, is once differentiable with continuous derivative on $[0,R+r]$, and moreover, has an infinitely differentiable and strictly concave restriction to $[T,R+r]$.    \textit{These additional properties may hold in greater generality.}

\end{document}